\documentclass[10pt,a4paper,twoside,leqno]{article}
\usepackage{amsfonts,amssymb}
\usepackage{eufrak}
\usepackage{amscd}

\newtheorem{defn}{Definition}

\vspace{0.5cm}

 4
\font\ebf=cmbx8
\font\erm=cmr8

\usepackage{graphicx}

\begin{document}

\thispagestyle{empty}

\noindent {\bf More on Combinatorial  interpretation of the
fibonomial coefficients }

\vspace{0.7cm} {\it Andrzej K. Kwa\'sniewski}

\vspace{0.7cm}

{\erm Higher School of Mathematics and Applied Informatics}

\noindent{\erm  Kamienna 17, PL-15-021 Bia\l ystok , Poland}

\vspace{0.5cm}

\noindent{\erm  ArXiv :  math.CO/0403017  v 1  1 March 2004}

\noindent {\ebf Summary}

{\small Combinatorial interpretation of the  fibonomial
coefficients recently attampted by the  present  author [1,2] and
presented here with suitable improvements results in a proposal of
a might be combinatorial interpretation of the recurrence relation
for fibonomial coefficients . The presentation is provided within
the context of the classical combinatorics` attitude to that type
of basic enumerative problems [3,4]. This note apart from plane
grid coordinate system used is fitted with several figures (Fig.1
- corrected) and examples  which illustrate the exposition of
statements and an interpretation of the recurrence itself.

\vspace{0.7cm}

\section{Introduction}
 There are various classical interpretations of
binomial coefficients,  Stirling numbers, and the $q$- Gaussian
coefficients . Recently Kovanlina  have discovered  [3,4] a
 simple and natural unified combinatorial interpretation  of all
 of them in terms of object  selection from weighted boxes with
 and without box repetition. So we are now in a position of the
 following recognition:

The classical, historically established standard interpretations might
be schematically presented for the sake of hint as follows:

SETS : Binomial coefficient $\Big({n \atop k}\Big )$
 , $\Big ({(n+k-1)\atop k}\Big ) $ denote number of subsets (without
  and with repetitions) - i.e. we are dealing with LATTICE of
  subsets.

SET PARTITIONS:  Stirling numbers of the second kind  $\Big\{ {n
\atop k}\Big\} $    denote number of partitions into exactly $k$
blocs - i.e. we are dealing with LATTICE of partitions.

PERMUTATION PARTITIONS : Stirling numbers of the first kind $\Big[
{n \atop k}\Big] $ denote number of permutations containing
exactly $k$ cycles

SPACES: $q$-Gaussian coefficient $\Big({n \atop k}\Big )_q $
denote number of $k$-dimensional subspaces in $n-th$ dimensional
space over Galois field  $GF(q)$ [5,6,7] i.e. we are dealing with
LATTICE of subspaces.

\vspace{3mm}

 \textbf{Before Konvalina combinatorial interpretation}.

\vspace{3mm}

 Algebraic similarities of the above classes of situations provided
 Rota [5] and Goldman and Rota [8,9] with an incentive to start the algebraic
unification that captures the intrinsic properties of these
numbers. The binomial coefficients, Stirling numbers and Gaussian
coefficients appear then as the coefficients in the characteristic
polynomials of geometrical lattices  [5] (see also [10] for the
subset-subspace analogy). The generalized coefficients [3]  are
called Whitney numbers of the first (characteristic polynomials)
and the second kind (rank polynomials).

\vspace{3mm}

 \textbf{Konvalina combinatorial interpretation}

 All these cases above and the case of Gaussian coefficients of the first
 kind $q^{\Big({n \atop 2}\Big )}\Big({n \atop k}\Big )_q $ are given unified
 Konvalina combinatorial interpretation in terms  of  the \textbf{generalized} binomial
 coefficients of the first and of the second kind  (see: [3,4] ).

\textbf{Unknowns ?} As for the distinguished [11,12,13,14]
Fibonomial coefficient defined below - no combinatorial
interpretation was  known till  today now to the present author .
The aim of this note is to promote a long  time waited for
 - classical in the spirit - combinatorial interpretation of Fibonomial coefficients.
 Namely we propose following [1,2] such a partial ordered set that the Fibonomial coefficients
  count the number of specific finite ``$birth-selfsimilar$''  sub-posets of an infinite
  locally finite not of binomial type , non-tree poset naturally related to
  the Fibonacci tree of rabbits growth process. This partial ordered set is defined equivalently
  via $\zeta$  characteristic matrix of partial order relation from its Hasse diagram.
  The classical scheme to be continued through  "Fibonomials" interpretation is the following:

POSET :   Fibonomial coefficient $\left( \begin{array}{c} n\\k\end{array}
\right)_{F}$  is the number of  ``birth-selfsimilar'' subposets.

\vspace{3mm}

\section{Combinatorial Interpretation}
It pays to get used to write $q$ or $\psi$ extensions of binomial symbols in
 mnemonic convenient  upside down notation [16,17] .
\begin{equation}\label{eq1}
\psi_n\equiv n_\psi , x_{\psi}\equiv \psi(x)\equiv\psi_x ,
 n_\psi!=n_\psi(n-1)_\psi!, n>0 ,
\end{equation}
\begin{equation}\label{eq2}
x_{\psi}^{\underline{k}}=x_{\psi}(x-1)_\psi(x-2)_{\psi}...(x-k+1)_{\psi}
\end{equation}
\begin{equation}\label{eq3}
x_{\psi}(x-1)_{\psi}...(x-k+1)_{\psi}= \psi(x)
\psi(x-1)...\psi(x-k-1) .
\end{equation}
You may consult [16,17] for further development and profit from
the use of this notation . So also here we use this upside down
convention for Fibonomial  coefficients:

\vspace{3mm}
$$
\left( \begin{array}{c} n\\k\end{array}
\right)_{F}=\frac{F_{n}!}{F_{k}!F_{n-k}!}\equiv
\frac{n_{F}^{\underline{k}}}{k_{F}!},\quad n_{F}\equiv F_{n}\neq 0, $$

\noindent where we make an analogy driven [16,17]  identifications $(n>0)$:
$$
n_{F}!\equiv n_{F}(n-1)_{F}(n-2)_{F}(n-3)_{F}\ldots 2_{F}1_{F};$$
$$0_{F}!=1;\quad n_{F}^{\underline{k}}=n_{F}(n-1)_{F}\ldots (n-k+1)_{F}. $$

\noindent This is the specification of the notation from [16] for the
purpose  Fibonomial Calculus case (see Example 2.1 in [17]).

\vspace{3mm}

    Let us now define the partially ordered infinite set  $P$.   We
shall label its vertices  by pairs of coordinates: ${\langle i , j
\rangle} \in {N \times N_0}$ where $N_0$ denotes the nonnegative
integers. Vertices show up in layers ("generations") of $N \times
N_0$ grid  along the recurrently emerging subsequent $s-th$ levels
$\Phi_s$ where  $s\in N_0$  i.e.

\begin{defn}
$$\Phi_s =\{\langle j, s\rangle 1\leq j \leq s_F\}, {s\in N_0}. $$
\end{defn}

We shall refer to $\Phi_s$  as to  the set of vertices at the
$s-th$ level. The population of the  $k-th$ level ("generation" )
counts  $k_F$  different member vertices for $k>0$ and one for
$k=0$.

\vspace{2mm}

Here down a disposal of vertices on $\Phi_k$ levels is visualized.

\vspace{5mm}

$---\Uparrow-----\Uparrow----up --Fibonacci---stairs--\star--k-th-level$\\

$---- and ----- so ---- on ---- up    --- \Uparrow ----------$\\
$\star \star \star \star \star \star \star \star \star \star \star \star \star \star \star \star \star \star\star \star \star \star \star \star \star \star \star \star \star \star \star \star \star \star \star \star \star \star \star \star \star \star \star \star --\star \star \star \star\star10-th-level$\\
$\star \star \star \star \star \star \star \star \star \star \star \star \star \star \star \star \star \star \star \star \star \star \star \star \star \star \star \star \star \star \star\star\star\star----------- 9-th-level$\\
$\star \star \star \star \star \star \star \star \star \star \star \star \star \star \star \star \star \star \star \star \star------------------8-th-level$\\
$\star \star \star \star \star \star \star \star \star \star \star \star \star -------------------------7-th-level$\\
$\star \star \star \star \star \star \star \star-----------------------------6-th-level$\\
$\star \star \star \star \star ---------------------------------5-th-level$\\
$\star \star \star ---------------------------------- 4-th-level$\\
$\star \star -----------------------------------3-rd-level  $ \\
$\star ------------------------------------ 2-nd-level$\\
$\star ----------------------------------- 1-st-level$\\
$\star ----------------------------------- 0-th-level$\\
 \vspace{2mm}
   \textbf{Figure 0. The $s-th$ levels in $ N\times N_0 $}
   , $N_0$ - nonnegative integers

 \vspace{2mm}

Accompanying the set $E$ of edges to the set $V$ of vertices - we
obtain the Hasse diagram where here down ${p,q,s}\in N_0 $.

Namely
\begin{defn}

$$P=\langle V,E\rangle ,\\ V=\bigcup_{0\leq p}\Phi_p ,\\ E
=\{\langle\langle j , p\rangle ,\langle q ,(p+1) \rangle
\rangle\}\bigcup\{\langle\langle 1 , 0\rangle ,\langle 1 ,1
\rangle \rangle\} , 1 \leq j \leq {p_F} , 1\leq q \leq
{(p+1)_F}.$$
\end{defn}
\begin{defn}
The prototype cobweb sub-poset is : $P_m = \bigcup_{0\leq s\leq
m}\Phi_s.$
\end{defn}

In reference [2] a partially ordered infinite set $P$ was
introduced  via  descriptive picture of  its Hasse diagram. Indeed
, we may picture out the partially ordered infinite set $P$ from
the  Definition $1$ with help of the sub-poset $P_{m}$ ({\it
rooted at $F_{0}$ level of the poset}) to be continued then ad
infinitum in now obvious way as seen from the  $Fig.1$ of $P_{5}$
below. It looks like the Fibonacci rabbits` tree with a specific
``cobweb''. \vspace{2mm}

\begin{center}

\includegraphics[width=75mm]{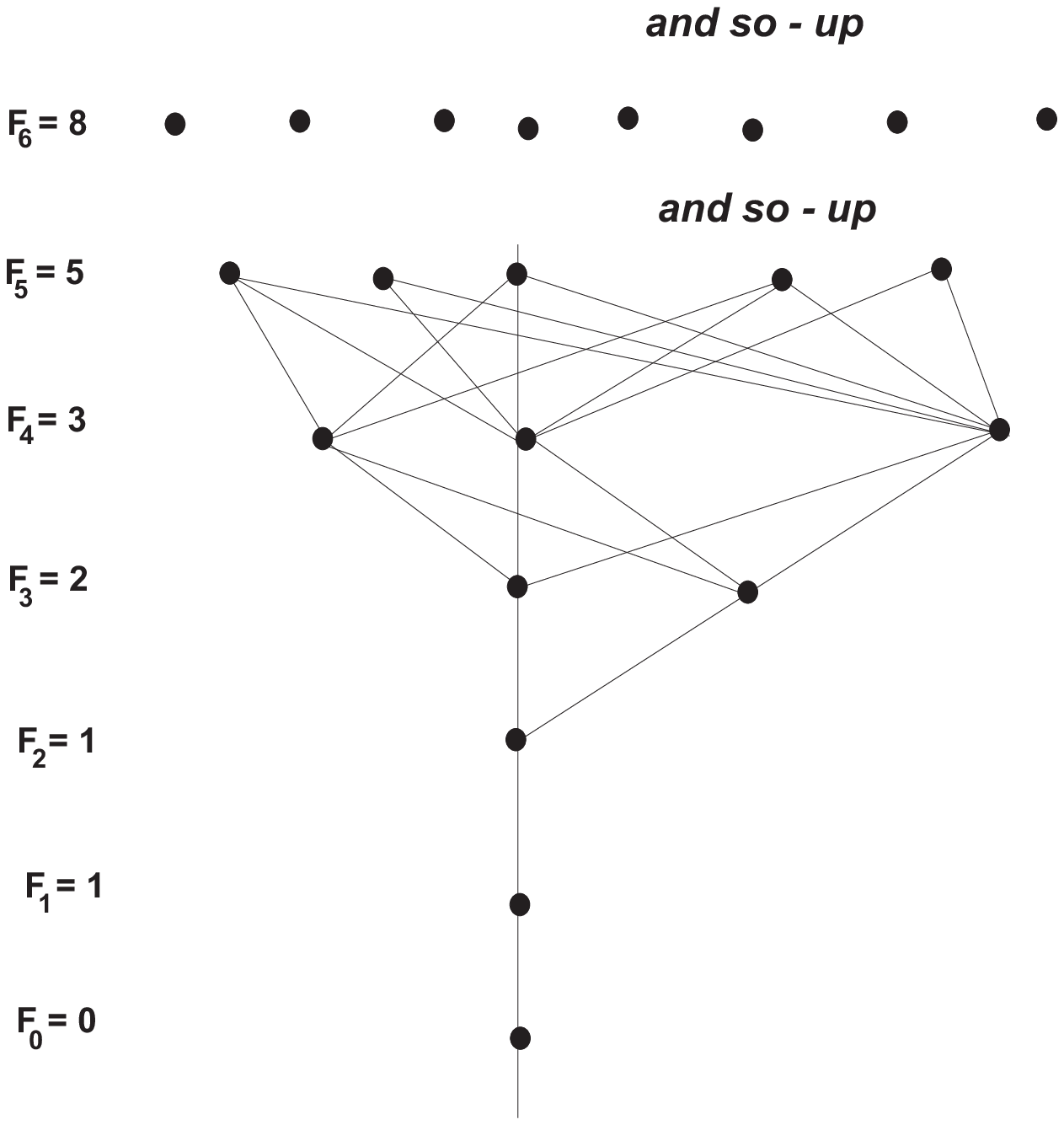}

\vspace{2mm}

\noindent {\small Fig.~1. Combinatorial interpretation of
Fibonomial coefficients.} \end{center}

\vspace{2mm} \noindent As seen above the $Fig.1$. displays the
rule of the construction of the  Fibonacci "cobweb"  poset. It is
being visualized clearly while defining this cobweb poset  $P$
with help of its incidence matrix . The incidence  $\zeta$
function [5,6,7] matrix representing uniquely just this cobweb
poset $P$ has the staircase structure correspondent with
$"cobwebed"$ Fibonacci Tree i.e. a Hasse diagram of the particular
partial order relation under consideration. This is seen below on
the Fig.$2$.

\vspace{2mm}

$$ \left[\begin{array}{ccccccccccccccccc}
1 & 1 & 1 & 1 & 1 & 1 & 1 & 1 & 1 & 1 & 1 & 1 & 1 & 1 & 1 & 1 & \cdots\\
0 & 1 & 1 & 1 & 1 & 1 & 1 & 1 & 1 & 1 & 1 & 1 & 1 & 1 & 1 & 1 & \cdots\\
0 & 0 & 1 & 1 & 1 & 1 & 1 & 1 & 1 & 1 & 1 & 1 & 1 & 1 & 1 & 1 & \cdots\\
0 & 0 & 0 & 1 & 0 & 1 & 1 & 1 & 1 & 1 & 1 & 1 & 1 & 1 & 1 & 1 & \cdots\\
0 & 0 & 0 & 0 & 1 & 1 & 1 & 1 & 1 & 1 & 1 & 1 & 1 & 1 & 1 & 1 & \cdots\\
0 & 0 & 0 & 0 & 0 & 1 & 0 & 0 & 1 & 1 & 1 & 1 & 1 & 1 & 1 & 1 & \cdots\\
0 & 0 & 0 & 0 & 0 & 0 & 1 & 0 & 1 & 1 & 1 & 1 & 1 & 1 & 1 & 1 & \cdots\\
0 & 0 & 0 & 0 & 0 & 0 & 0 & 1 & 1 & 1 & 1 & 1 & 1 & 1 & 1 & 1 & \cdots\\
0 & 0 & 0 & 0 & 0 & 0 & 0 & 0 & 1 & 0 & 0 & 0 & 0 & 1 & 1 & 1 & \cdots\\
0 & 0 & 0 & 0 & 0 & 0 & 0 & 0 & 0 & 1 & 0 & 0 & 0 & 1 & 1 & 1 & \cdots\\
0 & 0 & 0 & 0 & 0 & 0 & 0 & 0 & 0 & 0 & 1 & 0 & 0 & 0 & 1 & 1 & \cdots\\
0 & 0 & 0 & 0 & 0 & 0 & 0 & 0 & 0 & 0 & 0 & 1 & 0 & 1 & 1 & 1 & \cdots\\
0 & 0 & 0 & 0 & 0 & 0 & 0 & 0 & 0 & 0 & 0 & 0 & 1 & 1 & 1 & 1 & \cdots\\
0 & 0 & 0 & 0 & 0 & 0 & 0 & 0 & 0 & 0 & 0 & 0 & 0 & 1 & 0 & 0 & \cdots\\
0 & 0 & 0 & 0 & 0 & 0 & 0 & 0 & 0 & 0 & 0 & 0 & 0 & 0 & 1 & 0 & \cdots\\
0 & 0 & 0 & 0 & 0 & 0 & 0 & 0 & 0 & 0 & 0 & 0 & 0 & 0 & 0 & 1 & \cdots\\
. & . & . & . & . & . & . & . & . & . & . & . & . & . & . & . & . \cdots\\
 \end{array}\right]$$

 \vspace{2mm}

\textbf{Figure 2.  The  staircase structure  of  incidence matrix
$\zeta$ }

\vspace{2mm}

\textbf{Description} The "cob-viewer" encounters : in the zeroth -
zero and in the  $k-th$ row (for $k>0$) $F_k -1$ zeros right to
the diagonal value $1$ thus getting a picture of descending down
to infinity led by diagonal direction with use of growing in size
cobweb Fibonacci staircase tiled and build of the only up the
diagonal zeros
 - note these are forbiddance zeros (they code no edge links along $k-th$ levels
 ("generations") of the "$cobwebed$"  Fibonacci rabbits tree from [2].

\vspace{2mm}

This staircase structure of incidence [6,7] matrix $\zeta$ which
equivalently defines uniquely  this particular cobweb poset was
being recovered right from the Definition 1 and illustrative Hasse
diagram in Fig.$1$ of Fibonacci cobweb  poset. Let us say it again
- if one decides to define the poset $P$ by incidence matrix
$\zeta$  then must arrives at $\zeta$ with this easily
recognizable staircase-like structure of zeros in the upper part
of this upper triangle incidence matrix $\zeta$ just right from
the picture (see [1]  and [18] for recent references).

Let us recall [5,6,7] that $\zeta$ is being defined for any poset
as follows
 $(p,q \in P)$:
$$   \zeta (p,q)=\left\{ \begin{array}{cl} 1&for\ p \leq q,\\0&otherwise.
\end{array} \right. $$

\vspace{2mm}

The above $\zeta$ characteristic matrix of the partial order
relation in $P$ has been expressed explicitly in [1] in terms of
the infinite Kronecker delta  matrix $\delta$ from incidence
algebra $I(P)$ [5,6,7] as follows:

\vspace{2mm}
                   $$\zeta = \zeta_1 - \zeta_0 $$ where for $\langle x,y\rangle \in
{N_0\times N_0}$,
$$\zeta_1(x,y) = \sum_{k\geq 0} \delta (x+k,y)$$  while

$$\zeta_0(x,y)=\sum_{k\geq 0} \sum_{s\geq 1}
\delta (x,F_{s+1}+k)\sum_{1\leq r\leq (F_s
-k-1)}\delta(k+F_{s+1}+r,y).$$\\
 Naturally

$$   \delta(x,y)=\left\{ \begin{array}{cl} 1&for\ x = y,\\0&otherwise.

\end{array} \right. $$

\textbf{Important.} The knowledge of $\zeta$  matrix explicit form
enables one [6,7] to construct (count) via standard algorithms
[6,7]  the M{\"{o}}bius matrix $\mu =\zeta^{-1} $ and other
typical elements of incidence algebra perfectly suitable for
calculating number of chains, of maximal chains etc. in finite
sub-posets of $P$. Right from the definition of $P$ via its Hasse
diagram  here now follow quite obvious observations  .

\vspace{2mm}

\noindent {\bf Observation 1}

{\it The number of maximal chains starting from The Root  (level
$0_F$) to reach any point at the $n-th$ level  with $n_F$ vertices
is equal to $n_{F}!$}.

\vspace{2mm}

\noindent {\bf Observation 2} $(k>0)$

{\it The number of maximal chains \textbf{rooted in any fixed}
vertex at the $k-th$ level reaching the $n-th$ level with $n_F$
vertices is equal to $n_{F}^{\underline{m}}$, where $m+k=n.$ }

\vspace{2mm}

Indeed. Denote the number of ways to get along
maximal chains from a \textbf{fixed point} in $\Phi_k $to$  \Rightarrow  \Phi_n , n>k$ with the symbol\\
  $$[_{fixed}\Phi_k \rightarrow \Phi_n]$$
  then obviously we have :\\
           $$[\Phi_0 \rightarrow \Phi_n]= n_F!$$ and
$$[\Phi_0 \rightarrow \Phi_k]\times [_{fixed}\Phi_k\rightarrow \Phi_n]=
[\Phi_0 \rightarrow \Phi_n].$$

\vspace{2mm}

\textbf{Note} that the number  $[\Phi_k \rightarrow \Phi_n]$ of
all maximal chains starting from the $k-th$ level and ending at
the $n-th$ level equals to

$$ [\Phi_k \rightarrow \Phi_n]= k_F\times n_{F}^{\underline{m}}$$, where $m+k=n.$ \vspace{2mm}

   In order to find out the combinatorial interpretation of
Fibonomial coefficients let us consider all such  finite
\textit{"max-disjoint"} sub-posets rooted at the $k-th$ level at
any fixed vertex $\langle r,k \rangle, 1 \leq r \leq k_F $  and
ending  at corresponding number of vertices   at the $n-th$ level
($n=k+m$).

\vspace{2mm}

\textbf{Explanation:} \textit{"max-disjoint"} means : sub-posets
looked upon as families of \textbf{maximal} chains are disjoint
copies . These copies are like $P_m(k)_r$ copies isomorphic to
$P_m = P_m(0)_0$ - defined below for the sake of illustration.

In coordinate system this auxiliary illustrative cobweb sub-poset
$P_m(k)_r$ is defined as follows:

\begin{defn}
Let $\langle k,r\rangle \oplus \Pi$ denotes the shift of the set
$\Pi$ with the vector $\langle k,r\rangle$ Let $\Phi_o = \{\langle
0,0 \rangle\}.$
 Then we define:\\
$P_m(k)_r=\langle V_m(k)_r,  E_m(k)_r\rangle , V_m(k)_r= \langle
k,r\rangle \oplus \bigcup_{0\leq s\leq m}\Phi_s ,\\
E_m(k)_r =\{\langle\langle {(r+j)},(k+s)\rangle,\langle (r+i)
,(k+s+1)\rangle \rangle , 0\leq (r+j) \leq (k+s)_F , 1\leq {(r+i)}
\leq {(k+s+1)_F}\}.$
\end{defn}
\textbf{Observe}    $P_m(0)_0 = \langle V_m(0)_0, E_m(0)_0 \rangle
\equiv\langle V_m, E_m\rangle \equiv P_m .$ Hence $V_m(k)_r
=\langle k,r\rangle \oplus V_m $. Here, let us recall: $P_m$ is
the sub-poset of $P$ rooted at the $0-th$ level consisting of all
intermediate level vertices up to $m-th$ level ones - those from
$\Phi_m$ included $(See: Fig.1.)$.

\vspace{2mm}

\textbf{A newly k-th level born sub-cob browsing.}\\
 Consider now the following behavior of a sub-cob useful animal
moving from any given point of the $F_k$ "generation level"  of
the poset up and then up... It behaves as it has been born right
there and can reach at first $F_2$ vertices-points up, then $F_3$
points up , $F_4$ up... and so on - thus climbing up to the level
$F_{k+m} = F_n$ of the poset  $P$.  It can see  and then
potentially follow-  one of its own thus  accessible  isomorphic
\textit{copy} of sub-poset $P_m(k)_r$ in between the $k$-th and
$n$-th levels. ( "It" behaves exactly  as its Great Ancestor does
born at the Source Root $F_0-th$ level).

\vspace{2mm}

One of many of such \textit{"max-disjoint"} copies isomorphic
with the  sub-poset $P_m$`s (- the copies rooted at any fixed
point of the $k-th$ level) might be then found as a choice  to
start maximal chains forwarding up to the $n-th$ level - in the
limits of the chosen copy.

\vspace{2mm}

\textbf{How many} \textit{different} of such
\textit{"max-disjoint"} subposets \textit{choices} can be made?

\vspace{5mm}

\noindent {\bf Observation 3} $(\textbf{k}>\textbf{1})$

{\it Let  $n = k+m$. The number \textit{max-disjoint} sub-posets
isomorphic to $P_{m}$ , rooted  at the $k-th$ level and ending at
the n-th level  is equal to}
$$(n-m)_F\times\frac{n_{F}^{\underline{m}}}{m_{F}!} =
(n-m)_F\times\left( \begin{array}{c} n\\m\end{array}\right)_{F}$$
$$ =  k_F\times\left(\begin{array}{c} n\\k\end{array} \right)_{F}=
k_F\times\frac{n_{F}^{\underline{k}}}{k_{F}!}. $$

\vspace{2mm}

Indeed. Consider the number of all max-disjoint isomorphic copies
of $P_m$ rooted at a fixed vertex $\langle (r+j),k \rangle , 1\leq
(r+j) \leq k_F $. Denote this number with the symbol

$$ \left( \begin{array}{c} n\\k\end{array}\right)_{F}$$

Recall that  the number of all maximal chains from any point in
$\Phi_k $to$  \Rightarrow  \Phi_n , n>k$ is equal to\\
  $$ [\Phi_k \rightarrow \Phi_n]= k_F\times n_{F}^{\underline{m}}$$.

Then  one observes that :

\begin{equation}
k_F\times\left( \begin{array}{c} n\\k\end{array}\right)_{F} \times
[\Phi_0 \rightarrow \Phi_m] = [\Phi_k \rightarrow \Phi_n]=
k_F\times n_{F}^{\underline{m}}
\end{equation}

where  $[\Phi_0 \rightarrow \Phi_m]= m_F!$ counts the number of
maximal chains in any copy of the $P_m$. The factor  $k_F$ arises
from symmetric input by the vertices of the $k-th$ level.

\vspace{2mm}

For example: for the case $k=3 $ and $n=4$, according to the
interpretation given above, there should be $F_3\times {4 \choose
1}_F=F_3\times F_4! /F_1!F_3! = 6 $ max-disjoint copies of $P_1$
rooted  at the third level and  ending at the fourth level , which
is exact.

\vspace{2mm}

For example: for the  case $k=2$ and $n=4$, according to the
interpretation given above, there should be $F_2 \times {4 \choose
2}_F=F_4! /F_2!F_2! = 6$ max-disjoint copies of $P_2$ rooted at
the second level - ending at the fourth level , which is exact.

\vspace {2mm}

For example: for the  case $k=3$ and $n=5$, according to the
interpretation given above, there should be $F_3\times {5 \choose
3}_F = F_3 \times F_5! /F_3!F_2! = 30$ max-disjoint copies of
$P_2$ rooted at the third level and ending at the fifth level ,
which is exact.

\vspace{2mm}

For example: for the  case $k=2$ and $n=5$, according to the
interpretation given above, there should be $F_2\times {5 \choose
3}_F = F_5! /F_3!F_2! = 15$ max-disjoint copies of $P_3$ rooted at
the second level and ending at the fifth level , which is exact.

\vspace{2mm}

For example: for the  case $k=4$ and $n=5$, according to the
interpretation given above, there should be $F_4\times {5 \choose
1}_F =F_4 F_5! /F_4!F_1! = 15$ the number of max-disjoint copies
of $P_1$ rooted at the fourth level and ending at the fifth level
, which is exact.

\vspace{2mm}

\textbf{Important to Note} The reason for the restriction $k>1$
being applied is because $F_1=F_2$  and for $k=1$ Observation 3 is
not true.

\vspace{2mm}

For example: for the case $\textbf{k=1}$ and $n=4$, according to
the interpretation given above, there should be $F_1\times {4
\choose 1}_F= F_4! /F_1!F_3! = 3 $ max-disjoint copies of $P_3$
rooted at the thirst level and ending at the fourth level , which
is not true.

\vspace{2mm}

For example: for the  case $\textbf{k=1}$ and $n=5$, according to
the interpretation given above, there should be $F_1\times {5
\choose 1}_F = F_5! /F_4!F_1! = 5$ max-disjoint copies of $P_4$
rooted at the first level and ending at the fifth level which is
not true.

\renewcommand{\thesubsection}{\arabic{subsection}.}

\vspace{3mm}

\section{Does Konvalina like interpretation  of objects $F$- selections from weighted boxed exist?}

Binomial enumeration or finite operator calculus of Roman-Rota and Others is now the
standard tool of combinatorial analysis. The corresponding $q$-binomial calculus ($q$-calculus
- for short) is also the basis of much numerous  applications (see [19,20] for
altogether couple of thousands of respective references via enumeration and links).
In this context Konvalina unified binomial coefficients look intriguing and much promising.
The idea of  $F$-binomial or Fibonomial finite operator calculus (see Example 2.1 in [17])
consists of specification of the general scheme  - (see: [16,17] and references also to
Ward, Steffensen ,Viskov , Markowsky and others - therein)- specification  via the choice
 of the Fibonacci sequence to be sequence defining the generalized binomiality of polynomial
 bases involved (see Example 2.1 in [17]).Till now however we had been  lacking alike combinatorial
 interpretation of Fibonomial coefficients. We hope that this note would help not only via
 Observations above but also due to coming next- observation where recurrence relation
 for Finonomial coefficients is  is subjected to accordingly attempted combinatorial  interpretation.

\vspace{2mm}

\noindent{\bf Observation 4} $(k>0)$  ,  (combinatorial
interpretation of the recurrence)

The following known [11,14] recurrences hold
$$\left( \begin{array}{c} {n+1}\\k\end{array}
\right)_{F}=  F_{k-1} \left( \begin{array}{c} n\\k\end{array}
\right)_{F} + F_{n-k+2} \left( \begin{array}{c} n\\{k-1}\end{array}
\right)_{F}$$
or equivalently
$$\left( \begin{array}{c} {n+1}\\k\end{array}
\right)_{F}=  F_{k+1} \left( \begin{array}{c} n\\k\end{array}
\right)_{F} + F_{n-k} \left( \begin{array}{c} n\\{k-1}\end{array}
\right)_{F}$$ where
$$\left( \begin{array}{c} n\\0\end{array}
\right)_{F}= 1 , \left( \begin{array}{c} 0\\k\end{array}
\right)_{F}=0 , $$ \vspace{2mm} due to the recognition that we are
dealing with two disjoint classes in $P_{(n+1)}$ (n =k+m).
\textbf{The first one} for which
 $$ F_{k+1} \left( \begin{array}{c} n\\k\end{array}
\right)_{F}$$ equals to  $F_{k+1}$  times number of different
isomorphic copies of $P_m$ - rooted at a fixed point on the $k-th$
level (see Interpretation below) and

 \vspace{2mm}

\textbf{the second one} for which

\begin{equation}\label{eq4}
 F_{n-k} \left( \begin{array}{c} n\\{k-1}\end{array}\right)_{F}=
 F_{n-k} \left( \begin{array}{c} n\\{n-k+1}\end{array}\right)_{F}
\end{equation}

equals to $F_{n-k}$ times number of different isomorphic to
$P_m`s$ copies- rooted at a fixed point at the $(k-1)-th$ level
and ending at the $n-th$ level - see Interpretation below.

  \vspace{5mm}

 \textbf{Interpretation}   $(k>0)$

 \vspace{5mm}

$\diamond \diamond \diamond \diamond \diamond \diamond \diamond
\diamond \diamond \diamond F_{(n-k+2)} \diamond \diamond \diamond
\diamond \star\star\star\star\star\star\star \star-- \star
F_{(n+1)}-F_{(n-k+2)} \star\star\star\star\star\star\star
$\\

$\diamond\diamond\diamond\diamond\diamond\diamond\diamond\diamond\diamond\diamond\diamond\diamond\diamond\diamond\diamond\diamond\diamond\diamond\diamond\diamond\diamond---\diamond\diamond-----\star \star \star \star----n-th--level \star --- \star\star $\\

$\diamond\diamond\diamond\diamond\diamond\diamond\diamond\diamond\diamond\diamond\diamond\diamond\diamond\diamond\diamond\diamond\diamond\diamond\diamond\diamond\diamond\star\star\star\star\star\star\star\star\star\star\star\star\star\star\star\star\star\star\star\star\star\star\star--\star\star\star\star\star\star\star--\star\star\star\star\star\star\star$\\

$\diamond \diamond \diamond \diamond \diamond \diamond \diamond \diamond \diamond \diamond \diamond \diamond \diamond \star \star \star \star \star\star \star \star \star \star \star \star \star \star \star \star \star \star \star \star \star \star \star \star \star \star \star \star \star \star \star \star \star \star \star --10-th-level$ \\

$\diamond \diamond \diamond \diamond \diamond \diamond \diamond \diamond \star \star \star \star \star \star \star \star \star \star \star \star \star \star \star \star \star \star \star \star \star \star \star\star\star\star-------- 9-th-level$\\

$\diamond \diamond \diamond \diamond \diamond \star \star \star \star \star \star \star \star \star \star \star \star \star\star\star\star\star\star----------------8-th-level$\\

$\diamond\diamond\diamond \star \star \star \star \star \star \star \star \star \star -------------------------7-th-level$\\

$\diamond\diamond \star \star \star \star \star \star-----------------------------6-th-level$\\

$\diamond \star \star \star \star ---------------------------------5-th-level$\\

$\diamond \star \star ---------------------------------- 4-th-level$\\

$\star \star -----------------------------------3-rd-level  $ \\

$\star ------------------------------------ 2-nd-level$\\

$\star ----------------------------------- 1-st-level$\\

$\star --------------------------------- 0-th-level$\\

 \vspace{2mm}

   \textbf{Figure 3. The Diamond choice - two disjoint classes in $P_{(m+1)}(k)_r$.}

 \vspace{5mm}
The Fig.3 illustrates how the two disjoint classes referred to in
Observation 4  come into existence ($r=1, k=4$). First: every
cobweb sub-poset has the "trunk" of length $\geq$ one (in the
$Fig.3$ it is the extreme left maximal chain).  From any selected
root-vertex in $k-th$ level $F_{(k+1)}$ trunks may be continued in
$F_{k+1}$ ways. A trunk of the $P_{m+1}$ copy  being chosen - for
example the set of vertices $\langle1,s \rangle , k \leq s \leq
{(n+1)}$) in the case of diamond cobweb poset selected in Fig. 3 -
the resulting sub-cobweb ends with $F_{(n-k+2)}$ diamond vertices
("leafs") at $(n+1)-th$ level. The different copies when shifted
(in $F_{k+1}$ ways - each ) up and correspondingly completed by -
with the ultimate rightist maximal chain ending - the lacking part
of now $P_{m+1}`s$ copy - become now the different  copies rooted
at $k-th$ level and ending at the $(n+1)th$ level. This gives

 \vspace{5mm}

 $$ F_{k+1}\left(\begin{array}{c}
n\\k\end{array}\right)_{F} ,$$

 \vspace{2mm}

what constitutes  the first summand of the corresponding
recurrence.

 \vspace{5mm}

Consider then the non-$\Phi_k$ level (then to be shifted) choice
of the vertex .

 \vspace{10mm}

$\diamond \diamond \diamond \diamond \diamond \diamond \diamond
\diamond \diamond \diamond F_{(n+1)}-F_{(n-k+3)} \diamond \diamond
\diamond \diamond \star\star\star\star\star\star\star --\star
\otimes F_{(n-k+3)}
\otimes\otimes\otimes---\otimes\otimes\otimes\otimes\otimes$\\

$\diamond\diamond\diamond\diamond\diamond\diamond\diamond\diamond\diamond\diamond\diamond\diamond\diamond\diamond\diamond\diamond\diamond\diamond\diamond\diamond---\diamond\diamond-----\star \star \star -- \otimes--n-th--level \otimes --- \otimes\otimes $\\

$\diamond\diamond\diamond\diamond\diamond\diamond\diamond\diamond\diamond\diamond\diamond\diamond\diamond\diamond\diamond\diamond\diamond\diamond\diamond\diamond\diamond\star\star\star\star\star\star\star\star\star\star\star\star\star\star\star\star\star--\otimes\otimes\otimes\otimes\otimes\otimes\otimes--\otimes\otimes\otimes\otimes\otimes\otimes\otimes$\\

$\diamond \diamond \diamond \diamond \diamond \diamond \diamond \diamond \diamond \diamond \diamond \diamond \diamond \star \star \star \star \star\star \star \star -- \star \star \otimes\otimes\otimes\otimes\otimes\otimes\otimes\otimes\otimes\otimes\otimes\otimes\otimes\otimes\otimes\otimes\otimes\otimes\otimes\otimes\otimes 10-th-level$ \\

$\diamond \diamond \diamond \diamond \diamond \diamond \diamond \diamond \star \star \star \star \star \star \star  \star  \star  \star \star \otimes\otimes\otimes\otimes\otimes\otimes\otimes\otimes\otimes\otimes\otimes\otimes\otimes----------- 9-th-level$\\

$\diamond \diamond \diamond \diamond \diamond \star \star  \star  \star  \star \star \star \star \otimes\otimes\otimes\otimes\otimes\otimes\otimes\otimes------------------8-th-level$\\

$\diamond\diamond\diamond \star \star \star  \star  \star \otimes\otimes\otimes\otimes\otimes -----------------------7-th-level$\\

$\diamond\diamond \star \star \star \otimes\otimes\otimes----------------------------6-th-level$\\

$\diamond \star \star \otimes\otimes ---------------------------------5-th-level$\\

$\diamond \star \otimes---------------------------------- 4-th-level$\\

$\star \otimes--------------------------------3-rd-level$\\

$\otimes-------------------------------- 2-nd-level$\\

$\star ----------------------------------- 1-st-level$\\

$\star ----------------------------------- 0-th-level$\\

 \vspace{2mm}

   \textbf{Figure 4. The non-diamond choice - two disjoint classes in $P_{(n+1)}$.}

 \vspace{2mm}

\vspace{2mm}

The Fig.4 continues to illustrate how the two disjoint classes
referred to in Observation 4  are introduced. Now - what we do we
choose a vertex-root $\otimes$ in $(k-1)-th$ level in one of
$F_{(k-1)}$ ways. A trunk being chosen - say of the $\otimes$
cobweb sub-poset in Fig. 4 - it ends with $F_{(n-k+2)}$ $\otimes$
vertices ("leafs") at $n-th$ level. The number of $\otimes$
different isomorphic copies of $P_m$ rooted at a fixed point at
$(k-1)-th level$ is equal to  $ \left(\begin{array}{c}
n\\k-1\end{array}\right)_{F}$. These  copies  become now copies
rooted at $k-th$ level while ending at the $(n+1)-th$ level when
shifted up the $k-th$ and rooted at the same "diamond" root of the
first choice being afterwards  correspondingly completed by - with
the ultimate leftist maximal chain ending -the lacking part of
$P_{m+1}`s$  copies with all other copies rooted there at $k-th$
level and ending at the $(n+1)-th$ level. The number of thus
obtained different  copies is equal to

$$F_{n-k} \left( \begin{array}{c} n\\{k-1}\end{array}
\right)_{F}.$$

All together this gives the number of all different cobweb
sub-posets isomorphic copies  ending at $\Phi_{(n+1)}$ while
starting from a fixed point of  $\Phi_k$ level. This number  is
equal to the sum of number from the two disjoint classes i.e.

$$\left( \begin{array}{c} {n+1}\\k\end{array}
\right)_{F}=  F_{k+1} \left( \begin{array}{c} n\\k\end{array}
\right)_{F} + F_{n-k} \left( \begin{array}{c} n\\{k-1}\end{array}
\right)_{F}.$$

\vspace{5mm}

 In this connection the intriguing  question arises : May
one extend-apply somehow Konvalina theorem [3,4] below so as to
encompass also Fibonomial case under investigation ?

In [3,4] Konvalina considers  $n$ distinct boxes labeled with $i\in {[n]} , [n] \equiv \{1,...n\}$
such that each of $i-th$ box contains  $w_i$ distinct objects. John Konvalina uses the convention
$1\leq w_1 \leq w_2\leq ...\leq w_n$. Vector $N^n \ni \vec w =(w_1, w_2 ,...,w_n)$ is the weigh
 vector then. Along with Konvalina  considerations we have from [3]:

\vspace{3mm} \textbf{The Konvalina Theorem 1}

\vspace{2mm}

Let  $\vec w =(w_1, w_2 ,...,w_n)$ where $ 1\leq w_1 \leq w_2 \leq ...\leq w_n$.  Then

\textbf{I.} $$  C_{k}^n(\vec w) = C_{k}^{n-1}(\vec w) + w_n C_{k-1}^{n-1}(\vec w)$$

\textbf{II.}$$  S_{k}^n(\vec w) = S_{k}^{n-1}(\vec w) + w_n S_{k-1}^{n-1}(\vec w).$$

Here  $C_{k}^n(\vec w)$  denotes the generalized binomial coefficient of the first
kind with weight $\vec w$ and  it is the number of ways to select $k$ objects from
$k$ (necessarily distinct !) of the $n$ boxes with constrains as follows :  choose
$k$ distinct labeled boxes $$ i_1<i_2 <...<i_k $$  and then choose one object
 from each of the $k$ distinct boxes selected.  Naturally  one then has [3]

$$C_{k}^n(\vec w) = \sum_{1\leq i_1< i_2 < ...< i_k\leq n} w_{i_1}w_{i_2}...w_{i_k} . $$

Complementarily   $S_{k}^n(\vec w)$ denotes the generalized binomial coefficient of the second
kind with weight $\vec w$ and  it is the number of ways to select $k$ objects from
$k$ (not necessarily  distinct)  of the $n$ boxes with constrains as follows [3]:  choose $k$
not necessarily distinct labeled boxes $$ i_1\leq i_2 \leq ...\leq i_k$$  and then choose
one object from  each of the $k$ (not necessarily  distinct) boxes selected.  Obviously
one then has

$$S_{k}^n(\vec w) = \sum_{1\leq i_1\leq i_2 \leq ...\leq i_k\leq n} w_{i_1}w_{i_2}...w_{i_k}.$$

here the natural question arises : how are we to extend  Konvalina
theorem [3,4] above so as to encompass also Fibonomial case under
investigation ?

\vspace{2mm}

\textbf{Information I }: \textit about the preprint [21] entitled
\textit{Determinants, Paths,  and Plane Partitions} by Ira M.
Gessel,  X. G. Viennot [21].

Right after  Theorem $25$ - Section 10 , page 24 in [21])-
relating  the number $ N(R)$ of nonintersecting $k$-paths to
Fibonomial coefficients via $q$-weighted type counting formula-
the authors express their  wish worthy to be quoted: "\textit{it
would be nice to have a more natural interpretation then the one
we have given}"... " \textit{ R. Stanley has asked if there is a
binomial poset associated with the Fibonomial coefficients}..." -
Well. The cobweb locally finite infinite poset  by Kwasniewski
from [15,18,1,2] is not of binomial type. Recent incidence algebra
origin arguments [22] seem to make us not to expect  binomial type
poset come into the game.
    The $q$-weighted type counting formula from [21]gives rise to
an interesting definition of Fibonomial coefficients all together
with its interpretation in terms of nonintersecting $k$-paths due
to the properties of binomial determinant. Namely ,following [21]
let us consider points  $P_i= \langle\ 0,-i \rangle$ and $Q_i =
\langle\-n+i,-n+i \rangle$. Let  $R=\{r_1 < r_2 < ...<r_k\}\equiv
R(\vec r)$ be a subset of $\{0,1,...,n\}\equiv[n+1]$. Let $N(R)$
denotes the number of non-intersecting $k$-paths from $\langle
P_{r_1},...,P_{r_k}\rangle$ to $\langle
Q_{r_1},...,Q_{r_k}\rangle$ . Then $det \left( \begin{array}{c}
{r_i}\\n-r_{k+1-j}\end{array} \right) $ = $N(R)$.  The
$q$-weighted type counting formula from [21]then for $q=1$ means
that

$$\left( \begin{array}{c} {n+1}\\k\end{array}
\right)_{F}= \sum_{R(\vec r)}N(R) .$$

    In view of [21] another question arises - what is the relation like between
these two: Gessel and Viennot [21] non-intersecting $k$-paths and
cobweb sub-poset [18,2,3] points of view?

\vspace{2mm}

\textbf{Information II }: \textit on the partial ordered poset and
Fibonacci numbers paper [23] by Istvan Beck. The author of [23]
shows that  $F_n$ equals to the number of of ideals in a simple
poset called "fence" . This allows Him to infer via combinatorial
reasoning the identities :

$$F(n) = F(k) F(n + 1 - k) + F(k - 1)F(n -k)$$
$$F(n)= F(k-1) F(n + 1 -( k-1)) + F(k - 2) F(n -(k-1)).$$
 A straightforward application
of these above is the confirmation - just by checking - the
intriguing validity of recurrence relation for Fibonomial
coefficients . As we perhaps might learn from this note coming to
the end - both the Fibonomial coefficients as well as their
recurrence relation  are interpretable along the classical
historically established manner referring to the number of
objects` choices - this time these are partially ordered sub-sets
here called the cobweb sub-posets - the effect of the diligent
spider`s spinning of the maximal chains cobweb   during the
arduous day spent on the infinite Fibonacci rabbits` growth tree.

 \vspace{2mm}

\textbf{{\bf} Historical Memoir Remark}
The Fibonacci sequence origin is attributed and referred to the first
edition (lost)  of  ``Liber abaci'' (1202) by Leonardo Fibonacci  [Pisano]
(see second edition from 1228 reproduced as Il Liber Abaci di Leonardo
Pisano publicato secondo la lezione Codice Maglibeciano by  Baldassarre
Boncompagni  in Scritti di Leonardo Pisano  vol. 1, (1857) Rome).

\vspace{2mm}

\textbf{{\bf} Historical Quotation Remark } As accurately noticed
by Knuth and Wilf in [14]  the recurrent relations for Fibonomial
coefficients appeared already in $1878$  Lukas work [11]. In our
opinion - Lucas`s Th\'eorie des fonctions num\'eriques simplement
p\'eriodiques is the far more non-accidental context for binomial
and binomial-type coefficients - Fibonomial coefficients included.

While studying this mentioned important and inspiring paper by
Knuth and Wilf [14] and in the connection with a context of this
note a question raised by the authors with respect to their
formula (15) is worthy to be repeated : \textit{Is there a
"natural" interpretation....} - May be then fences from [23]  or
cobweb posets or  ... "Natural" naturally might have many
effective faces ...

\vspace{2mm}

 \textbf{ Acknowledgements} I am very much indebted to Mgr Ewa Krot -  for
 her substantial critical remarks allowing to present my exposition
 in hopefully better shape.

 \vspace{2mm}

\begin
{thebibliography}{99}
\parskip 0pt

\bibitem{1}
A. K. Kwa\'sniewski, {\it The logarithmic Fib-binomial formula}
Advanced Stud. Contemp. Math. {\bf 9} No 1 (2004):19-26

\bibitem{2}
A. K. Kwasniewski {\it Information on combinatorial interpretation
of Fibonomial coefficients }   Bull. Soc. Sci. Lett. Lodz Ser.
Rech. Deform. 53, Ser. Rech.Deform. {\bf 42} (2003): 39-41 ArXiv:
math.CO/0402291   v1 18 Feb 2004

\bibitem{3}
J. Konvalina , {\it Generalized binomial coefficients and the
subset-subspace  problem } , Adv. in Appl. Math. {\bf 21}  (1998)
: 228-240

\bibitem{4}
J. Konvalina , {\it A Unified Interpretation of the Binomial
Coefficients, the Stirling Numbers and the Gaussian Coefficients}
The American Mathematical Monthly {\bf 107}(2000):901-910

\bibitem{5}
 Gian-Carlo Rota "On the Foundations of Combinatorial Theory,
I. Theory  of  Möbius Functions";   Z. Wahrscheinlichkeitstheorie
und Verw. Gebiete, vol.2 ,  1964 , pp.340-368.

\bibitem{6}
 E. Spiegel, Ch. J. O`Donnell  {\it Incidence algebras}  Marcel
Dekker, Inc. Basel $1997$ .

\bibitem{7}
Richard P. Stanley {\it Enumerative Combinatorics}  {\bf I},
Wadsworth and Brooks Cole Advanced Books  and Software, Monterey
California, $1986$

\bibitem{8}
 Goldman J.  Rota G-C. {\it The Number of Subspaces in a vector
space} in  Recent Progress in Combinatorics (W. Tutte, Ed.):
75-83, Academic Press, San Diego, $1969$, see: ( J. Kung, Ed.)
"Gian Carlo Rota on Combinatorics"  Birkhäuser, Boston
(1995):217-225

\bibitem{9}
Goldman J.  Rota G-C. {\it On the Foundations of Combinatorial
Theory IV; finite-vector spaces and Eulerian generating functions}
Studies in Appl. Math. {\bf 49} (1970): 239-258

\bibitem{10}
J. Kung   {\it The subset-subspace analogy} ( J. Kung, Ed.)
"Gian Carlo Rota on Combinatorics"  Birkhäuser, Boston (1995):277-283

\bibitem{11}
E. Lucas, {\it Th\'eorie des fonctions num\'eriques simplement
p\'eriodiques}, American Journal of Mathematics {\bf 1} (1878):
184--240; (Translated from the French by Sidney Kravitz), Ed. D.
Lind, Fibonacci Association, 1969.

\bibitem{12}
 G. Fonten\'e, {\it G\'en\'eralisation d`une formule connue},
Nouvelles Annales de Math\'ematiques (4) {\bf 15} (1915), 112.

\bibitem{13}
H. W. Gould,   {\it The bracket function and Fonten\'e-Ward
generalized binomial coefficients with applications to Fibonomial
coefficients}, The Fibonacci Quarterly {\bf 7} (1969), 23--40.

\bibitem{14}
D. E. Knuth, H. S. Wilf  {\it The Power of a Prime that Divides a
Generalized Binomial Coefficient} J. Reine Angev. Math. {\bf 396}
(1989): 212-219

\bibitem{15}
A. K. Kwa\'sniewski, {\it More on Combinatorial interpretation of
Fibonomial coefficients},   Inst. Comp. Sci.  UwB/Preprint no. 56,
November 2003.

\bibitem{16}
A. K. Kwa\'sniewski,   {\it Towards  $\psi$-extension of finite
operator calculus of Rota}, Rep. Math. Phys. {\bf 47} no. 4
(2001), 305--342.    ArXiv: math.CO/0402078  2004

\bibitem{17}
A. K. Kwa\'sniewski, {\it On simple characterizations of Sheffer
$\Psi$-polynomials and related propositions of the calculus of sequences},
Bull.  Soc.  Sci.  Lettres  \L \'od\'z {\bf 52},S\'er. Rech. D\'eform.
 {\bf 36} (2002), 45--65. ArXiv: math.CO/0312397  $2003$

\bibitem{18}
A. K. Kwa\'sniewski,{\it Comments on  combinatorial interpretation
of fibonomial coefficients - an email  style letter} , Bulletin of
the Institute of Combinatorics and its Applications, {\bf 42}
September (2004): 10-11

\bibitem{19}
T. Ernst  , {\it The History of q-Calculus and a new Method },
 http://www.math.uu.se/~thomas/Lics.pdf 19 December (2001),
(Licentiate Thesis). U. U. D. M. Report (2000).

\bibitem{20}
A. K. Kwa\'sniewski,{\it First Contact Remarks on Umbra Difference Calculus References Streams}
 Inst. Comp. Sci.  UwB Preprint No  {\bf 63} (January $2004$ )

\bibitem{21}
Ira M. Gessel,  X. G. Viennot{\it Determinant Paths and Plane
Partitions  } preprint  (1992)
http://citeseer.nj.nec.com/gessel89determinants.html

\bibitem{22}
A.K.Kwasniewski {\it The second part of on duality triads`
paper-On fibonomial and other triangles  versus  duality triads}
 Bull. Soc. Sci. Lett. Lodz Ser. Rech. Deform. 53, Ser. Rech.
Deform. {\bf 42} (2003): 27 -37  ArXiv: math.GM/0402288 v1 18 Feb.
$2004$

\bibitem{23}
I. Beck  {\it Partial Orders and the Fibonacci Numbers} The
Fibonacci Quarterly {\bf 26} (1990): 172-174 .

\end{thebibliography}



\end{document}